\documentclass[a4paper]{amsart}
\usepackage{amsmath,amsfonts,amssymb,amsthm}
\usepackage{latexsym, xspace, enumerate}
\usepackage[mathscr]{eucal}
\usepackage[all]{xy}
\usepackage{hyperref}
\usepackage{amscd}
\usepackage{mathtools}
\usepackage{stmaryrd}
\usepackage{rotating}
\usepackage{xfrac}

\newtheorem{theorem}{Theorem}[section]

\newtheorem{corollary}[theorem]{Corollary}

\theoremstyle{definition}

\newtheorem{question}[theorem]{Question}

\newtheorem*{theoremKK}{Kaloujnine-Krasner Theorem}
\newtheorem*{theoremE}{Extension Theorem}
\newtheorem{thmx}{Theorem}

\newcommand{\C}{\mathbb C}

\def\Sym{\mathrm{Sym}}

\newcommand{\Sp}{\text{Sp}}

\newcommand{\wrwr}{\operatorname{\wr\wr}}

\makeatletter
\def\Ddots{\mathinner{\mkern1mu\raise\p@
\vbox{\kern7\p@\hbox{.}}\mkern2mu
\raise4\p@\hbox{.}\mkern2mu\raise7\p@\hbox{.}\mkern1mu}}
\title{Unrestricted wreath products and sofic groups}
\author[G. Arzhantseva]{Goulnara Arzhantseva}
\address[Goulnara Arzhantseva]{Universit\"{a}t Wien, Fakult\"{a}t f\"{u}r Mathematik, Oskar-Morgenstern-Platz~1, 1090 Wien, Austria}
\email[Goulnara Arzhantseva]{goulnara.arzhantseva@univie.ac.at}

\author[F. Berlai]{Federico Berlai}
\address[Federico Berlai]{UPV/EHU, Department of Mathematics, Barrio Sarriena s/n, 48940, Leioa, Spain}
\email[Federico Berlai]{federico.berlai@ehu.eus}

\author[M. Finn-Sell]{Martin Finn-Sell}
\address[Martin Finn-Sell]{Universit\"{a}t Wien, Fakult\"{a}t f\"{u}r Mathematik, Oskar-Morgenstern-Platz~1, 1090 Wien, Austria}
\email[Martin Finn-Sell]{martin.finn-sell@univie.ac.at}

\author[L. Glebsky]{Lev Glebsky}
\address[Lev Glebsky]{Universidad Aut\'{o}noma de San Luis Potos\'{i}, Mexico}
\email[Lev Glebsky]{glebsky@cactus.iico.uaslp.mx}


\begin{document}
\begin{abstract}
We show that the unrestricted wreath product of a sofic group by an amenable group is sofic. We use this result
to present an alternative proof of the known fact that any group extension with sofic kernel and amenable quotient is again a sofic group. 
Our approach exploits the famous Kaloujnine-Krasner theorem and extends, with an additional argument,
to hyperlinear-by-amenable groups.
\end{abstract}
\maketitle
\section{Introduction}
We propose a new proof of the following fact, originally proved by Elek and Szab\'{o} \cite{ES}.
\begin{theoremE}\emph{
Let $G$ be a group with a sofic, normal subgroup $N\trianglelefteqslant  G$ such that the quotient $G/N$ is amenable. Then $G$ is sofic. }
\end{theoremE}
Our proof relies on considering unrestricted wreath products. Given two groups $G$ and $H$, 
their \emph{unrestricted wreath product} $G\wrwr H$ is, by definition, the semi-direct product
\begin{equation*}
G\wrwr H:=\Bigl(\prod_{h\in H}G\Bigr)\rtimes H,
\end{equation*}
where $H$ acts on the direct product shifting the coordinates:
\begin{equation*}
\bigl((g_x)_{x\in H},h\bigr) \cdot \bigl((g'_x)_{x\in H},h'\bigr)=\bigl((g_{h'x}g'_{x})_{x\in H},hh'\bigr). 
\end{equation*}
Subgroups of the unrestricted wreath product $G\wrwr H$ play a decisive role in understanding the extensions that can be obtained starting from the two groups $G$ and~$H$:
\begin{theoremKK}\emph{
Let $G$ and $H$ be groups, and let $\pi\colon G\wrwr H \twoheadrightarrow H$ denote the canonical projection onto the quotient $H$. There is a bijection of sets between
\begin{equation*}
\{[E]_{\cong} \mid 1 \to G\hookrightarrow E \twoheadrightarrow H\to 1\}
\end{equation*}
and
\begin{equation*}
\{[E]_{\text{conj}}\mid E\leqslant G\wrwr H,\,  \pi(E)=H,\, E\cap \ker \pi\cong G \}
\end{equation*}
that induces a group isomorphism between the extension $E$ and the subgroup $E$ of $G\wrwr H$.}
\end{theoremKK}
We prove the following theorem.
\begin{thmx}\label{thm:unrestricted}\emph{
Let $G$ be a sofic group and $H$ be an amenable group. Then the unrestricted wreath product $G\wrwr H$ is sofic.}
\end{thmx}

This allows us to conclude the Extension Theorem, independently of \cite{ES}, by
combining the Kaloujnine-Krasner Theorem and the obvious fact that soficity is preserved under taking subgroups.

We adapt our strategy, using an additional argument, to hyperlinear groups.

\begin{thmx}\label{thm:hyper}\emph{
Let $G$ be a hyperlinear group and $H$ be an amenable group. Then the unrestricted wreath product $G\wrwr H$ is hyperlinear.}
\end{thmx}

As hyperlinearity is preserved under taking subgroups, using the Kaloujnine-Krasner Theorem, we immediately obtain 

\begin{corollary}
Let $G$ be a group with a hyperlinear, normal subgroup $N\trianglelefteqslant  G$ such that the quotient $G/N$ is amenable. Then $G$ is hyperlinear. 
\end{corollary}
 This fact is perhaps known to specialists but it had not yet appeared in the literature explicitly.

In the final section we indicate that the results are optimal for such an approach: it does not immediately extend neither 
to sofic-by-$\{$residually amenable$\}$ nor to $\{$weakly sofic$\}$-by-amenable group extensions.

\subsection*{Acknowledgements}
This work greatly benefited from the \textquotedblleft Measured group theory\textquotedblright\ program held in 2016 and 
the  \textquotedblleft Word maps and stability of representations\textquotedblright\
 workshop held in 2013 at the Erwin Schr\"{o}dinger International Institute for Mathematics and Physics, Vienna. We thank the institution
for its warm hospitality and for the excellent working environment. 

Goulnara Arzhantseva and Martin Finn-Sell are partially supported, Federico Berlai, and Lev Glebsky's visit to Vienna in 2016, were supported by the European Research Council (ERC) grant of Prof. Goulnara Arzhantseva, grant agreement no.~259527. Lev Glebsky's visit to Vienna in 2013
was supported by the Austrian-Swiss research network on \textquotedblleft Sofic groups\textquotedblright\, grant CRSI22-130435 of the Swiss NSF.
Federico Berlai is now supported in part by ERC grant PCG-336983, Basque Government Grant IT974-16, and Ministry of Economy, Industry and Competitiveness of the Spanish Government Grant MTM2017-86802-P.
\section{Sofic groups}\label{section1}

Sofic groups were introduced in 1999 by Gromov \cite{Gro} (under a different terminology, as groups with \emph{initially subamenable Cayley graphs}\footnote{One should not be confused with Gromov's, more restrictive, definition of \emph{initially subamenable groups}.}) in his work on Gottschalk's surjunctivity conjecture.
The terminology ``sofic''  is due to Weiss \cite{Weiss} and comes from the Hebrew word ``finite''. Sofic groups are a common generalization of
residually finite and amenable groups. They form 
a wide 
subclass of so-called \emph{hyperlinear} groups which appeared in operator algebra in the context of Connes' embedding conjecture. Both sofic and hyperlinear groups are subject of current intensive research, see~\cite{Pes08,CSC,CL} and references therin. 
  
To define sofic groups, we recall the notion of normalised Hamming distance on a finite symmetric group. 
Let $F$ be a finite set and $\lvert F\rvert$ denote its cardinality. The following bi-invariant distance is called the normalised \emph{Hamming distance} on $\Sym(F)$: for permutations $\alpha, \beta\in \Sym(F),$ we define
\begin{equation*}\label{definition_Hamming}
d_F(\alpha,\beta):=\frac{1}{\lvert F\rvert}\Bigl\lvert\Bigl\{f\in F\mid \alpha(f)\neq\beta(f)\Bigr\}\Bigr\rvert.
\end{equation*}
A group $G$ is said to be \emph{sofic} if for every finite subset $K\subseteq G$ and for every $\varepsilon>0$ there exist a finite set $F$ and a map 
$\varphi\colon G\to \Sym(F)$ satisfying:
\begin{enumerate}
 \item[(s1)] for all $k_1,k_2\in K$ we have that $d_{F}\bigl(\varphi(k_1k_2),\varphi(k_1)\varphi(k_2)\bigr)\leqslant\varepsilon$;
 \item[(s2)] for all distinct $k_1,k_2\in K$ we have that $d_F\bigl(\varphi(k_1),\varphi(k_2)\bigr)\geqslant 1-\varepsilon$.
\end{enumerate}
We call such a map $\varphi\colon G\to\Sym(F)$ a $(K,\varepsilon)$-\emph{approximation}, and when $K$ and $\varepsilon$ are clear from the context, we 
refer to it simply as to an \emph{approximation}.
Notice that, instead of condition $(s2)$, one could ask for $d_F\bigl(\varphi(k_1),\varphi(k_2)\bigr)$ to be bounded away from $0$ by a constant $\alpha$ independent from the finite set $K$ (confront~\cite[Proposition~3.4]{GR}).

\medskip
Amenable groups are particular instances of sofic groups. Indeed, let $K$ be a finite subset of an amenable group $G$ and consider $\varepsilon>0$. 
Then, as $G$ is amenable, there exists a \emph{F\o{}lner set} corresponding to $K$, i.e. a
non-empty finite subset $F\subseteq G$ such that
\begin{equation}\label{amenable_eq1}
\lvert F\setminus sF\rvert<\varepsilon'\lvert F\rvert,\qquad \forall s\in S:=\bigl(\{e_G\}\cup K\cup K^{-1}\bigl)^2, 
\end{equation}
where $\varepsilon':=\varepsilon /\lvert S\rvert$. Without loss of generality, it may be supposed that the F\o{}lner set $F$ is symmetric, i.e. $F=F^{-1}$.

As shown, for example, in \cite[Proposition 7.5.6]{CSC}, in this case, there exist a finite set $E$ such that
\begin{equation}
E\subseteq F\subseteq G,\qquad \text{}\qquad \lvert E\rvert\geqslant(1-\varepsilon)\lvert F\rvert, 
\end{equation}
and a map $\varphi\colon G\to\Sym(F)$, defined by
\begin{equation}\label{sofic.approximation.general}
\varphi(g)f:=\begin{cases}
               gf&\text{if }gf\in F\\
               \alpha_g(gf)&\text{if }gf\in gF\setminus F,
              \end{cases}
\end{equation}
where $\alpha_g\colon gF\setminus F\to F\setminus gF$ is an arbitrarily chosen fixed bijection, such that
\begin{equation}\label{eq_AmenableCondition1}
\bigl(\varphi(k_1)\varphi(k_2)\bigr)f=k_1k_2f=\varphi(k_1k_2)f,\qquad \forall k_1,k_2\in K,\ \forall f\in E\subseteq F,
\end{equation}
and
\begin{equation}\label{eq_AmenableCondition2}
\varphi(k_1)f=k_1f\neq k_2f=\varphi(k_2)f,\qquad \forall k_1,k_2\in K, k_1\neq k_2,\ \forall f\in E\subseteq F,
\end{equation}
and so that $\varphi$ is a $(K, \varepsilon)$-approximation.

We stress that the preceding two conditions are quite strong because $F$ is a subset of $G$: we see from Equations \eqref{eq_AmenableCondition1} and \eqref{eq_AmenableCondition2} that
the approximation $\varphi$ is governed by the group operation on a big part of the set $F$. In general, $F$ might not be a subset of the group we want to approximate, and it might not be equipped with a partial operation induced by the group. 

\section{Proof of Theorem~\ref{thm:unrestricted}}
We now prove that, given a sofic group $G$ and an amenable group $H$, the unrestricted wreath product $G\wrwr H$ is sofic.

Consider a finite subset $K\subseteq G\wrwr H$ and fix $\varepsilon>0$. Without loss of generality, we assume that the identity element of $G\wrwr H$ belongs to $K$, and that $0<\varepsilon<1$.

We define the maps
\[\pi_1\colon G\wrwr H \twoheadrightarrow \prod_HG,\qquad \pi_2\colon G\wrwr H\twoheadrightarrow H\]
to be the projections on the first and second component of the semi-direct product, respectively. Note that the map $\pi_1$ is not a group homomorphism.

As the set $K$ is finite, there exists a finite separating set $I\subseteq H$ for $\pi_1(K)\subseteq \prod_HG$. This means that, given two distinct sequences
$(g_x)_{x\in H}$ and $(g'_x)_{x\in H}$ in $\pi_1(K)$, there exists an index $i\in I$ such that $g_i\neq g'_i$.

Since the group $H$ is amenable, its approximations can be constructed starting from appropriate F\o{}lner sets as discussed in the previous section.

Let
\begin{equation*}    
K_H:=\pi_2(K)\cup I
\end{equation*}
and let us choose an $\varepsilon_H>0$ such that 
\begin{equation}\label{varepsilon_H}
\varepsilon_H<\varepsilon. 
\end{equation}
As recalled in Section \ref{section1}, there exist a finite, symmetric subset $B\subseteq H$ and a $(K_H,\varepsilon_H)$-approximation $\varphi_H\colon H\to \Sym(B)$, such that
\begin{equation}\label{sofic.approximation.H}
\varphi_H(h)b=\begin{cases}
               hb&\text{if }hb\in B\\
               \alpha_h(hb)&\text{if }hb\in hB\setminus B,
              \end{cases}
\end{equation}
where $\alpha_h\colon hB\setminus B\to B\setminus hB$ is a fixed bijection. Moreover, there exists a set $E_B\subseteq B$ such that $\lvert E_B\rvert\geqslant(1-\varepsilon_H)\lvert B\rvert$ and, 
for all $h,h'\in K_H$ and for all $b\in E_B$, we have that
\begin{equation}\label{sofic.approximation.H.explicit}
\varphi_H(h)b=hb,\qquad \bigl(\varphi_H(h)\varphi_H(h')\bigr)b=hh'b=\varphi_H(hh')b. 
\end{equation}
By Equation \eqref{amenable_eq1}, as $B$ is a F\o{}lner set corresponding to $K_H$, we have that
\begin{equation}\label{equation.folner}
\lvert B\setminus hB\rvert< \varepsilon_H\lvert B\rvert,\qquad \forall h\in K_H=\pi_2(K)\cup I.
\end{equation}
As concerns the approximation of $G$, let us consider
\begin{equation*}
K_G:=\{g_x\in G\mid (g_x)_{x\in H}\in\pi_1(K),\, x\in B^{-1}\cdot B\}. 
\end{equation*}
Choose an $\varepsilon_G>0$ such that
\begin{equation}\label{varepsilon_G.second.condition}
1-(1-\varepsilon_H)(1-\varepsilon_G)^{\lvert B\rvert}\leqslant \varepsilon.
\end{equation}
This condition is equivalent to require that
\begin{equation}\label{eq_tricky_sofic}
(1-\varepsilon_G)^{\lvert B\rvert}\geqslant\frac{1-\varepsilon}{1-\varepsilon_H}.
\end{equation}
As $\varepsilon_H<\varepsilon$ by Equation \eqref{varepsilon_H}, the right hand side of Equation \eqref{eq_tricky_sofic} is strictly smaller than one.
Hence, as the constants $\varepsilon$ and $\varepsilon_H$ - and consequently also $\lvert B\rvert$ - are fixed, there exists $\varepsilon_G>0$ satisfying Equation~\eqref{varepsilon_G.second.condition}.

Since $G$ is sofic, there exist a finite set $A$ and a $(K_G,\varepsilon_G)$-approximation $\varphi_G\colon G\to\Sym(A)$. Moreover, for each choice of $g,g'\in K_G$ there exists 
a finite set $E_A=E_A^{g,g'}\subseteq A$, which depends on $g$ and $g'$, such that $\lvert E_A\rvert\geqslant(1-\varepsilon_G)\lvert A\rvert$ and
\begin{equation}\label{sofic.approximation.G.explicit}
\bigl(\varphi_G(g)\varphi_G(g')\bigr)a=\varphi_G(gg')a,\qquad \forall a\in E_A. 
\end{equation}
Given the approximations $\varphi_G$ and $\varphi_H$ we construct an approximation $\varphi$ for $G \wrwr H$. Consider the finite set $C:=B\times A^B,$ where $A^B$ denotes the set of all functions from $B$ to $A.$
We define the map $\varphi\colon G\wrwr H\to\Sym(C)$ as follows. For an element $(g,h)\in G\wrwr H$, with $g=(g_x)_{x\in H}\in \prod_HG$ and $h\in H$, and $(b,\tau)\in C$, let
\begin{equation}\label{definition.sofic.approximation.wreath}
\varphi(g,h)(b,\tau):=\bigl(\varphi_H(h)b,\bar\tau\bigr), 
\end{equation}
where $\bar\tau\colon B\to A$ is the function defined by
\begin{equation}\label{tau.bar.definition}
\bar\tau(i):=\varphi_G(g_{bi})\tau(i),\qquad \forall i\in B.
\end{equation}
We now prove that $\varphi$ satisfies condition (s1)  in the definition of soficity. We have to show that for $(g,h)$ and $(g',h')$ in the finite set $K\subseteq G\wrwr H$ the inequality
\begin{equation*}
d_C\bigl(\varphi(g,h)\varphi(g',h'),\varphi((g_{h'x}g'_{x}),hh')\bigr)\leqslant \varepsilon
\end{equation*}
holds.
Indeed, let $(b,\tau)\in C$. We have that
\begin{equation}
\varphi(\underbrace{(g_{h'x}g'_{x})}_{=(p_x)},hh')(b,\tau)=(\varphi_H(hh')b,\tau_1), 
\end{equation}
where
\begin{equation}\label{hom.tau1}
\tau_1(i)=\varphi_G(p_{bi})\tau(i)=\varphi_G(g_{h'bi}g'_{bi})\tau(i). 
\end{equation}
Moreover,
\begin{equation}
\bigl(\varphi(g,h)\varphi(g',h')\bigr)(b,\tau)=\varphi(g,h)\bigl(\varphi_H(h')b,\tau_2)\bigr)=\bigl((\varphi_H(h)\varphi_H(h'))b,\tau_3\bigr), 
\end{equation}
where
\begin{equation}
\tau_2(i)=\varphi_G(g'_{bi}) \tau(i)
\end{equation}
and
\begin{equation}\label{hom.tau3}
\tau_3(i)=\varphi_G(g_{(\varphi_H(h')b)i})\tau_2(i)=\varphi_G(g_{(\varphi_H(h')b)i})\bigl(\varphi_G(g'_{bi}) \tau(i)\bigr)
\end{equation}
Notice that if $b\in E_B$ (see Equation \eqref{sofic.approximation.H.explicit}), then $\varphi_H(h)b=hb\in B$. In this case, from Equation~\eqref{hom.tau3} we obtain that
\begin{equation}\label{hom.tau3emezzo}
\tau_3(i)= \varphi_G(g_{(\varphi_H(h')b)i})\bigl(\varphi_G(g'_{bi}) \tau(i)\bigr)=\bigl(\varphi_G(g_{h'bi})\varphi_G(g'_{bi})\bigr) \tau(i).
\end{equation}
Confronting Equations \eqref{hom.tau1} and \eqref{hom.tau3emezzo}, we notice that the functions $\tau_1$ and $\tau_3$ will be equal if, in addition to 
$b\in E_B$, we have that $\tau(i)\in E_A$ for all $i\in B$, see Equation~\eqref{sofic.approximation.G.explicit}.

Hence, we conclude that
\begin{equation*}\label{equation_with_troubles}
\begin{split}
d_C\bigl(\varphi(g,h)\varphi(g',h'),\varphi((g_{h'x}g'_{x}),hh')\bigr)&\leqslant 1- \frac{1}{\lvert C\rvert}\Bigl(\lvert E_B\rvert \cdot \lvert E_A\rvert^{\lvert B\rvert}\Bigr)\\ 
&=1-\frac{\lvert E_B\rvert}{\lvert B\rvert}\cdot\Bigl(\frac{\lvert E_A\rvert}{\lvert A\rvert}\Bigr)^{\lvert B\rvert}\\
&\leqslant 1-(1-\varepsilon_H)(1-\varepsilon_G)^{\lvert B\rvert}.
\end{split}
\end{equation*}
This is smaller or equal to $\varepsilon$ in view of Equation \eqref{varepsilon_G.second.condition}.

\smallskip
It remains to check that $\varphi$ also satisfies condition (s2) in the definition of soficity. 
Let $(g,h)$ and $(g',h')$ be distinct elements of the finite set $K\subseteq G\wrwr H$, 
with $g=(g_x)_{x\in H}$ and  $g'=(g'_x)_{x\in H}$. Two cases may occur: $h\neq h'$ or $h=h'$.

If $h\neq h',$ then
\begin{equation}\label{case.equal}
d_C\bigl(\varphi(g,h),\varphi(g',h')\bigr)\geqslant d_B\bigl(\varphi_H(h),\varphi_H(h')\bigr)\geqslant 1-\varepsilon_H> 1-\varepsilon,
\end{equation}
where the last inequality is due to Equation \eqref{varepsilon_H}.

If $h=h',$ then $(g_x)\neq (g'_x)$, and hence there exists an index $i$ in the finite set $I$ for which $g_i\neq g'_i$. 
To proceed, for almost all choices of $b\in B$, we want to express this coordinate $i\in I\subseteq H$ as $bb_i$, where $b_i\in B$ is an element that depends on $i$ and on the chosen $b$. We have that
\begin{equation}\label{I_need_B_symmetric}
i=bb_i\quad \Leftrightarrow\quad b=ib_i^{-1}\quad \Leftrightarrow\quad b\in B\cap iB^{-1}=B\cap iB,
\end{equation}
where in the equality of Equation \eqref{I_need_B_symmetric} we used the fact that the set $B$ is symmetric, $B^{-1}=B$.
Using Equation \eqref{equation.folner}, we obtain
\begin{equation*}
\lvert B\cap iB\rvert=\lvert B\rvert -\lvert B\setminus iB\rvert\geqslant \lvert B\rvert-\varepsilon_H\lvert B\rvert=(1-\varepsilon_H)\lvert B\rvert.
\end{equation*}
This means that, given $i\in I$, this coordinate can be expressed as $i=bb_i$ with $b,b_i\in B$ for almost all choices of $b\in B$ (`almost all' states for $\geqslant (1-\varepsilon_H)\lvert B\rvert$ choices of $b\in B$).

Suppose now that $(b,\tau)\in C$ and $b$ is such that $i=bb_i$, for $b_i\in B$. We have that
\begin{equation*}
\varphi(g,h)(b,\tau)=(\varphi_H(h)b,\tau_4),
\end{equation*}
where $\tau_4(j)=\varphi_G(g_{bj})\tau(j)$ for all $j\in B$, and
\begin{equation*}
\varphi(g',h)(b,\tau)=(\varphi_H(h)b,\tau_5),
\end{equation*}
where $\tau_5(j)=\varphi_G(g'_{bj})\tau(j)$ for all $j\in B$.

Since $G$ is sofic, $i=bb_i$ and $g_{bb_i},g'_{bb_i}$ are distinct elements, there exists $X\subseteq A$ such that $\lvert X\rvert\geqslant(1-\varepsilon_G)\lvert A\rvert$ for which
\begin{equation*}
\varphi_G(g_{bb_i})(a)\neq \varphi_G(g'_{bb_i})(a),\qquad \forall a\in X.
\end{equation*}
Therefore, the two functions $\tau_4$ and $\tau_5$ are distinct whenever $\tau(b_i)\in X$.
This implies that
\begin{equation}\label{case.different}
\begin{split}
d_C\bigl(\varphi(g,h),\varphi(g',h)\bigr)&\geqslant \frac{1}{\lvert C\rvert}\Bigl((1-\varepsilon_H)\lvert B\rvert \cdot (1-\varepsilon_G)\lvert A\rvert \cdot \lvert A\rvert^{\lvert B\rvert -1}\Bigr)\\
&=\frac{\lvert B\rvert\cdot\lvert A\rvert^{\lvert B\rvert}(1-\varepsilon_H)(1-\varepsilon_G)}{\lvert B\rvert\cdot\lvert A\rvert^{\lvert B\rvert}}\\
&=(1-\varepsilon_H)(1-\varepsilon_G)\geqslant 1-\varepsilon,
\end{split}
\end{equation}
where the last inequality is due to Equation \eqref{varepsilon_G.second.condition} and to the fact that $1-\varepsilon_G$ is bigger than $(1-\varepsilon_G)^{\lvert B\rvert}$.

To conclude, combine Equation \eqref{case.equal} and Equation \eqref{case.different} to see that $\varphi$ satisfies condition $(s2)$ in the definition of soficity, as required.

\section{Hyperlinear groups and concluding remarks}
In this final section we analyse our approach to metric approximations of extensions through those of unrestricted wreath products
in more general settings. Namely, we discuss two renown generalisations of sofic groups:  \emph{hyperlinear} groups and  \emph{weakly sofic}  groups.
Then we consider sofic-by-\{residually amenable\} group extensions and  conclude with remarks and questions on abelian-by-sofic extensions.

\subsection{Hyperlinear groups and proof of Theorem~\ref{thm:hyper}}
A very natural question is whether or not our argument extends to hyperlinear groups. 

A group $G$ is said to be \emph{hyperlinear} if for every finite subset $K\subseteq G$ and for every $\varepsilon>0$ there exist a natural number $n$ and a map $\vartheta\colon G\to \mathcal{U}(\C^n)$ satisfying:
\begin{itemize}
 \item[$(h1)$] for all $k_1,k_2\in K$ we have that $d_{HS}\bigl(\vartheta(k_1)\vartheta(k_2),\vartheta(k_1k_2)\bigr)\leqslant \varepsilon$;
 \item[$(h2)$] for all distinct $k_1,k_2\in K$ we have that $d_{HS}\bigl(\vartheta(k_1),\vartheta(k_2)\bigr)\geqslant \sqrt{2}-\varepsilon$.
\end{itemize}
Here, $d_{HS}$ is the normalised \emph{Hilbert-Schmidt distance} on $\mathcal{U}(\C^n)$, defined to be, for  $u=(u_{i,j}),v=(v_{i,j})\in \mathcal{U}(\C^n)$,
\begin{equation*}
d_{HS}(u,v):=\sqrt{\frac{1}{n}\sum_{i,j=1}^n \lvert u_{i,j}-v_{i,j}\rvert^2}.
\end{equation*}
Sofic groups are hyperlinear \cite{EShyp}. Indeed, it can be easily seen that
\begin{equation}\label{Hamming_HS}
d_F(\sigma,\tau)=\frac{1}{2}d_{HS}(P_\sigma,P_\tau)^2, 
\end{equation}
where $P_\sigma$ and $P_\tau$ are the permutation matrices in $\mathcal{U}(\C^{\lvert F\rvert})$ associated to permutations $\sigma,\tau\in\Sym(F)$ (see, e.g. \cite{EShyp}). The converse implication, i.e. whether or not hyperlinear groups are sofic, is a well-known open problem.
We refer to \cite{Pes08,CL} and to the references therein for more information on hyperlinearity.

\smallskip
Modifying the proof of Theorem \ref{thm:unrestricted}, together with an additional argument  below, we obtain Theorem~\ref{thm:hyper} from
the Introduction and, hence, the hyperlinearity of hyperlinear-by-amenable groups extensions.

\smallskip
The modifications required to our construction above are mostly minor and involve adjusting the conditions on $\varepsilon$ in terms of $\varepsilon_{G}$ and $\varepsilon_{H}$  in light of the relationship between Hamming distance and Hilbert-Schmidt distance. 

Fix an $\varepsilon>0$, $K\subseteq G \wrwr H$.
For an amenable group $H$, proceed as in the text above, find a F\o{}lner set $B$ for $K_{H}$ and an $\varepsilon_H>0$ satisfying
\begin{equation*}
 \sqrt{2\varepsilon_H}<\varepsilon.
\end{equation*}
Let $\vartheta_H\colon H\to \mathcal{U}(\C^{\lvert B\rvert})$ be the representation into $\mathcal{U}(\C^{\lvert B\rvert})$ corresponding to the sofic 
approximation $\varphi_H$:
\begin{equation*}
\vartheta_H(h):=P_{\varphi_H(h)}\in\mathcal{U}(\C^{\lvert B\rvert}).
\end{equation*}
As $G$ is supposed to be hyperlinear (instead of sofic), there exist a natural number $n$ and a $(K_G,\varepsilon_G)$-approximation $\vartheta_{G}\colon G\to \mathcal{U}(\C^n)$, for a constant $\varepsilon_{G}$ satisfying
\begin{equation*}
\sqrt{\frac{\lvert B\rvert}{\lvert B\rvert+n^{\lvert B\rvert}}\varepsilon_H^2+\frac{n^{\lvert B\rvert}}{\lvert B\rvert+n^{\lvert B\rvert}}
\varepsilon_G^{2\lvert B\rvert}}<\varepsilon.
\end{equation*}
To define
\begin{equation*}
\vartheta\colon G \wrwr H \rightarrow \mathcal{U}\bigl(\mathbb{C}^{\lvert B\rvert }\oplus\bigoplus_{B}\mathbb{C}^{n} \bigr),
\end{equation*}
let $\{e_b\}_{b\in B}$ be a orthonormal basis for $\C^{\lvert B\rvert}$, and for $v\in\bigoplus_B\C^{n}$, express it as $v = \sum_{i\in B} v_{i}$, 
where each $v_{i}$ is the projection onto the $i$-th copy of $\mathbb{C}^{n}$ in the summand.

Consider hence
\begin{equation*}
\vartheta(g,h)(e_{b},v) = (\vartheta_{H}(h)e_{b}, \overline{v}),
\end{equation*} 
where $\overline{v}=\sum_{i\in B}\overline{v}_i$
is defined by orthogonally summing the components
\begin{equation*}
\overline{v}_{i}:= \vartheta_{G}(g_{b^{-1}i})v_{i}.
\end{equation*}
Extend now $\vartheta(g,h)$ linearly in the first component to obtain the action on an arbitrary element $(w,v) \in \C^{\lvert B\rvert}\oplus\bigoplus_{B}\mathbb{C}^{n}$.
This gives a well-defined map and our prior calculations pass through without significant changes so that we obtain Theorem~\ref{thm:hyper}.
\subsection{Weakly sofic groups}
We observe now that the assumption on the soficity of the kernel group $G$ in Theorem~\ref{thm:unrestricted} could be considered optimal.

Following \cite{GR}, a group is said to be \emph{weakly sofic} if there exists $\alpha>0$ such that for every finite subset $K\subseteq G$ and for every $\varepsilon>0$ there exist a finite group $F$ equipped with a bi-invariant metric $d$ and  a 
map $\varphi\colon G\to F$ satisfying:
\begin{itemize}
 \item[$(ws1)$] for all $k_1,k_2\in K$ we have that $d\bigl(\varphi(k_1k_2),\varphi(k_1)\varphi(k_2)\bigr)\leqslant \varepsilon$;
 \item[$(ws2)$] for all distinct $k_1,k_2\in K$ we have that $d\bigl(\varphi(k_1),\varphi(k_2)\bigr)\geqslant \alpha$.
\end{itemize}
Sofic groups are weakly sofic as the normalised Hamming distance is bi-invariant. Again, the converse implication is a well-known open problem.

\smallskip
The proof of Theorem~\ref{thm:unrestricted} strongly relies on soficity, i.e. on approximations by finite symmetric groups equipped with the normalised Hamming distance. 
It does not apply in the a priori more general context of weakly sofic groups. 
The reason is that the approximation for the group $G\wrwr H$ constructed in Equation \eqref{definition.sofic.approximation.wreath} produces, starting from the maps $\varphi_G\colon G\to \Sym(A)$
and $\varphi_H\colon H\to \Sym(B)$, the new map $\varphi\colon G\wrwr H\to \Sym\bigl(B\times A^B\bigr)$. This new map is not diagonal, that is, its image is not contained in the subgroup
\begin{equation*}
\Sym(B)\times \Sym(A^B)\leqslant \Sym\bigl(B\times A^B\bigr).
\end{equation*}
This makes the subsequent arguments incompatible with the setting of arbitrary finite groups equipped with bi-invariant normalised metrics, other than symmetric groups with Hamming distances - which confirms that soficity is the optimal condition to expect in this context.

\subsection{Residual amenability}
We demonstrate here that the assumption on the amenability of the quotient group $H$ in Theorem~\ref{thm:unrestricted} could also be seen as optimal.

Indeed, suppose instead that $H$ is residually amenable. 
Using \cite[Proposition~7.3.2]{CSC}, we obtain a set $L$ with $K_H\subseteq L\subseteq H$ and a binary operation $\diamond$ such that 
$(L,\diamond)$ is an amenable group, and such that $\diamond$ corresponds to
the group operation of $H$ when restricted to the finitely many elements of $K_H$: $hh'=h\diamond h'$ for all $h,h'\in K_H$.

This causes several issues in the above computations. In particular, when comparing Equations \eqref{hom.tau1} and \eqref{hom.tau3},
one needs to restrict the attention to $b\in E_B\cap K_H$ (and not
just to $b\in E_B$, as in the amenable case). For everything to hold, we would need to impose, instead of Equation \eqref{varepsilon_G.second.condition}, the following new condition on $\varepsilon_G$:
\begin{equation*}\label{equation_nojoy}
1-\frac{\lvert E_B\cap K_H\rvert}{\lvert B\rvert}(1-\varepsilon_G)^{\lvert B\rvert}\leqslant \varepsilon. 
\end{equation*}
However, this new condition is equivalent to
\begin{equation*}
(1-\varepsilon_G)^{\lvert B\rvert}\geqslant  \frac{1-\varepsilon}{\lvert E_B\cap K_H\rvert /\lvert B\rvert}.
\end{equation*}
Therefore, we would require the quantity $(1-\varepsilon_G)^{\lvert B\rvert}$ (which is smaller than one) to be bigger than 
$\frac{1-\varepsilon}{\lvert E_B\cap K_H\rvert /\lvert B\rvert}$. However, this last quantity is bigger than one,
as $\lvert K_H\rvert /\lvert B\rvert\ll1-\varepsilon$.
This obvious contradiction confirms that the amenability assumption on the group $H$ is the suitable condition in this approach (as well as in the original proof of \cite{ES}), and it cannot be easily relaxed to residual amenability.

\subsection{Abelian-by-sofic extensions}
In a sense, extensions we consider in this paper are opposite to amenable-by-sofic groups, and nothing is known about metric approximations of such groups. For instance, it follows from \cite[Corollary 3.8]{Pau}
(see also \cite[Corollary 2]{HaSa}, where the authors focus on \emph{restricted} regular wreath products)
that if $F$ is a finitely generated free group with normal subgroup $S$, and $F/S$ is sofic, then also the group $F/S'$ is sofic, where $S'$ denotes the derived
subgroup of~$S$.

\smallskip
Let us consider a finitely generated group $G=F/R$, where $F$ is a finitely generated free group and $R$ is its normal subgroup. 
Let $N$ be an abelian normal subgroup of $G$, and let $Q=G/N$ be a sofic quotient of $G$. By construction, 
there exists a normal subgroup $S\trianglelefteqslant F$ such that $R\leqslant S$ and $Q=F/S$. Moreover, it follows that $N$ is isomorphic to $S/R$, and that the derived subgroup $S'$ of $S$ is contained in $R$, 
because $S/R$ is an abelian group.

Therefore, we obtain a well-defined surjective homomorphism $\pi\colon F/S'\twoheadrightarrow F/R$, and the following exact sequences:
\begin{equation*}
\xymatrix{ 
&   &         &  1       &       &1 \\
& 1\ar[r] & S/R\,\ar@{^{(}->}[r] & F/R\ar@{->>}[r]\ar[u]    & F/S\ar[r]\ar[ur] & 1\\
&   &         &  F/S' \ar@{->>}[u]^{\pi}\ar@{->>}[ur]  &       & \\
1&  \frac{S/S'}{R/S'}\ar[uur]_\cong \ar[l] &S/S'  \ar@{->>}[l]_{\quad\kappa}\ar@{^{(}->}[ur]       &  R/S'\ar@{_{(}->}[l]\ar@{^{(}->}[u]       & 1\ar[l]      & \\
&  1\ar[ur] &         &  1   \ar[u]    &       & 
}
\end{equation*}
A na\"{i}ve way to conclude that the abelian-by-sofic extension $G=F/R$ is sofic, is to require for $S'$ to coincide with $R$. In that case, $\kappa$ would be an isomorphism, and therefore $N\cong S/R$ would
be isomorphic to the abelianisation $S/S'$ of the free group $S$, that is, $N\cong \mathbb Z^s$, where $s$ is the rank of $S$.
This condition clearly excludes several cases, in particular when $N$ has torsion, or when its rank does not coincide with $s$.

A specific instance when $N$ indeed has torsion appears in the following question:
\begin{question}\label{question}
Let $N\trianglelefteqslant G$ be such that $N$ is finite cyclic and $G/N$ is a sofic (respectively: hyperlinear) group. Is $G$ sofic (respectively: hyperlinear)? 
\end{question}
Even the case when $\lvert N\rvert =2$ seems to be open.
Moreover, a positive answer to Question \ref{question} implies that any group extension $G$ with finitely generated abelian kernel and sofic quotient is itself sofic.
Indeed, any such $G$ is the limit, in the space of marked groups, of \{finite cyclic\}-by-sofic groups (for a proof of this fact, confront \cite[Theorem 2.1]{B}).

Groups which are worth mentioning explicitly in this context are Deligne's central extensions \cite{pD}:
\begin{equation*}
1\to \mathbb Z \hookrightarrow \widetilde{\Sp_{2n}(\mathbb Z)} \twoheadrightarrow \Sp_{2n}(\mathbb Z)\to 1,
\end{equation*}
where $n\geqslant 2$ and $\widetilde{\Sp_{2n}(\mathbb Z)}$ is the preimage of the simplectic group $\Sp_{2n}(\mathbb Z)$ in the universal cover $\widetilde{\Sp_{2n}(\mathbb R)}$ of $\Sp_{2n}(\mathbb R)$. They are not known to 
be sofic, or hyperlinear.

These groups are not residually finite, have Kazhdan's property~(T), and are finitely presented. Therefore, they are not initially subamenable (that is, they are not locally embeddable into amenable groups).
They are the limits in the space of marked groups, as $m\to\infty$, of the non-initially subamenable groups $\widetilde{\Sp_{2n}(\mathbb Z)}/m\mathbb Z$, for each fixed $n\geqslant 2$ and $m\in\mathbb N$, which are central group extensions of the form
\begin{equation*}
1\to {\mathbb Z}/{m\mathbb Z}  \hookrightarrow \widetilde{\Sp_{2n}(\mathbb Z)}/m\mathbb Z \twoheadrightarrow \Sp_{2n}(\mathbb Z)\to 1.
\end{equation*}
These central extensions with finite cyclic kernel are isolated points in the space of marked groups (although in \cite[Section 5.8]{CGP} slightly different groups are considered, we observe that the argument exploited there also proves that the above extensions with finite center are isolated).
Thus, a positive answer to Question \ref{question} restricted to these central extensions would then both settle the question of soficity (respectively: hyperlinearity) of Deligne's groups and give new examples of finitely presented non-initially subamenable groups admitting such metric approximations.

\end{document}